\documentclass[12pt,onecolumn,twoside]{IEEEtran}
\usepackage{amsmath,amssymb,euscript,psfrag,latexsym,psfrag,graphicx,float,color}
\usepackage{amsmath,amssymb,euscript,psfrag,latexsym,psfrag,graphicx,float,color}

\newtheorem{thm}{Theorem}

\newtheorem{prop}[thm]{Proposition}

\newtheorem{example}[thm]{Example}
\newcommand{\wno}{{\rm wno}}
\newcommand{\cG}{{\cal G}}
\newcommand{\bD}{{\bf D}}
\color{black}
\newcommand {\cL}{{\cal L}}
\newcommand {\cH}{{\cal H}}
\newcommand {\bP}{{\bf P}}
\newcommand {\bH}{{\bf H}}
\newcommand {\bC}{{\bf C}}
\newcommand {\bI}{{\bf I}}
\def\spacingset#1{\def\baselinestretch{#1}\small\normalsize}
\spacingset{1}

\title{\LARGE \bf
Feedback Control and the Arrow of Time
}

\author{Tryphon T. Georgiou and Malcolm C. Smith% <-this % stops a space
\thanks{This work was partially supported by the  National Science Foundation.
T.T. Georgiou is with Department of Electrical and Computer
Engineering, University of Minnesota, Minneapolis, MN 55455; {tryphon@ece.umn.edu}. 
M.C. Smith is with the Department of Engineering,
University of Cambridge, Cambridge, CB2 1PZ, U.K.;
{mcs@eng.cam.ac.uk}}%
}

\begin{document}

\maketitle
\thispagestyle{empty}
\pagestyle{empty}

%%%%%%%%%%%%%%%%%%%%%%%%%%%%%%%%%%%%%%%%%%%%%%%%%%%%%%%%%%%%%%%%%%%%%%%%%%%%%%%%
\begin{abstract}
  The purpose of this paper is to highlight the central role that the
  time asymmetry of stability plays in feedback control.  We show that
  this provides a new perspective on the use of doubly-infinite or
  semi-infinite time axes for signal spaces in control theory. We then
  focus on the implication of this time asymmetry in modeling
  uncertainty, regulation and robust control.  We point out that
  modeling uncertainty and the ease of control depend critically on
  the direction of time. We also discuss the relationship of this
  control-based time-arrow with the well known arrows of time in
  physics.
\end{abstract}

%\spacingset{1.9}
%%%%%%%%%%%%%%%%%%%%%%%%%%%%%%%%%%%%%%%%%%%%%%%%%%%%%%%%%%%%%%%%%%%%%%%%%%%%%%%%
\section{Introduction}

The origin and implications of the ``arrow of time'' is one of the
deepest and least understood subjects of physics. The ``arrow''
is an intrinsic part of the world as we know it. Yet its emergence in
thermodynamics and cosmology, from physical laws which are apparently
impervious to it, remains a controversial subject \cite{price}. At first sight,
this subject may seem unconnected with the theory of feedback control.
However, starting from the very basic fact that our notion of
stability in the sense of Lyapunov is time-asymmetric, we argue that
the ``arrow of time'' does have important implications on modeling and
uncertainty, robustness of stability, as well as on the topology for
the study of the dynamics of feedback interconnections.

The circle of ideas that gave rise to this paper began in a short note
published by the authors thirteen years ago \cite{GS95}. There, it was
pointed out that the doubly-infinite time axis presents some ``intrinsic
difficulties'' for developing a suitable input-output systems
theory---difficulties that are not present in the semi-infinite time
axis setting.  These difficulties are not mere mathematical
technicalities.  Rather, they relate fundamentally to the consistency
of the theory of stabilizability across different frameworks.
Subsequently, a number of papers were written which shed light on the
problem \cite{makila1,makila5,makila7,jacob,jacob2,jacob3}.  The
present paper takes a fresh look and traces the origin of the
``puzzle'' to the arrow of Lyapunov stability, and then, explores the
relevance of this arrow to the topology of dynamical systems and
feedback theory.

The relationship of the modern theory of dynamical systems with classical physics and thermodynamics is a developing one.
A classical contribution by Nyquist and Johnson \cite{nyquist,johnson} is a derivation of the electromotive force due to thermal agitation in conductors. In \cite{brockett} the issue of irreversibility is treated from the point of view of stochastic control theory. More recently \cite{haddad} has sought to formalize classical thermodynamics in the mathematical language of modern dynamical systems (see also \cite{delvenne}). In \cite{mitter} information flow and entropy have been studied in the context of the Kalman filter. In \cite{sandberg} it is shown that a linear macroscopic dissipative system can be approximated by a linear lossless microscopic system over arbitrary long time intervals. Our point of view here is influenced by \cite{price} and is somewhat different to the above references in that our main goal is to highlight a time-asymmetry, point out its implications, and discuss its relationship to other well-known asymmetries.

The present paper begins by providing a new explanation of the issues raised
in \cite{GS95} with regard to an input-output theory for the
doubly-infinite time axis. In Section \ref{optimal control} we
introduce the time-conjugation operator and discuss the implications of the time-arrow in optimal control problems.
In Section \ref{section:robustcontrol} we analyse the effect of the
time-arrow on modelling uncertainty; we show that dynamical systems
which are close in the usual sense, that a common controller can
stabilise and give similar closed-loop responses for either, may not
be close when the time-arrow is reversed.  Then, in Section
\ref{section:robustcontrol2}, we further illuminate the inherent
time-asymmetry in our ability to control a dynamical system with two
specific examples. These can be thought of as examples of time
irreversible feedback phenomena (see Section \ref{geddanken}).
In Section \ref{physics} we briefly discuss the arrow of time in physics and its relation to the time-arrow of feedback stability. Finally, in Section \ref{delays} we consider feedback loops with small time delays and discuss the contrasting effects of delays and predictors and the connection with the arrow of time.

\section{Time-asymmetry and stability}\label{section:timeasymmetry}
\subsection{Input-output and Lyapunov stability}
We focus on finite-dimensional linear dynamical systems which, for the
most part, are assumed to be time-invariant.  The dimensions of input,
state and output (column) vectors, as well as the consistent sizes of
transformation matrices in state-space models, are suppressed for
notational simplicity.  The following result is basic and well-known, cf.\ \cite[p. 52-53]{jlwillems}, \cite[p. 82]{GreenLimebeer}.
\begin{prop}\label{prop1} {\sf Let $\bf P$ be a linear time-invariant finite-dimensional system which is controllable and observable and is specified by
\begin{eqnarray}\label{system}
\dot{x} &=& Ax+Bu,\\
y&=&Cx+Du, \label{system1}
\end{eqnarray}
with an initial condition $x(0)=0$.  Then $y\in \cL_2[0,\infty)$ for
all $u\in \cL_2[0,\infty)$ if and only if the matrix $A$ is Hurwitz.
Moreover, if this condition holds, $y$ is determined uniquely by
$\hat{y}(s)=(C(sI-A)^{-1}B+D)\hat{u}(s)$, where $\hat{\;}$ denotes the
Laplace transform.}
\end{prop}

Many variants and extensions of the result are familiar: signal spaces
with different norms can also be used; there is a
finite-gain property relating the
$\cL_2$-norms of $y$ and $u$; even with $x(0)\neq 0$ the main
equivalence in the proposition still holds.
Here we would like to
highlight the fact that the result establishes an equivalence between
stability defined in terms of the {\em forced response} and stability
defined in terms of the {\em free response}, i.e. an equivalence
between bounded-input/bounded-output (BIBO) stability and Lyapunov
stability for a system operating on the positive time-axis.
Asymptotic stability in the sense of Lyapunov is obviously a
time-asymmetric concept since convergence of the state vector is
required as $t$ tends to PLUS infinity, starting from an arbitrary
initial condition at $t=0$. In itself, BIBO stability does not appear
to have this asymmetry, yet it is implicit in the formulation of
Proposition~\ref{prop1}.

To further illustrate the point we can write down the following
obvious corollary of Proposition~\ref{prop1}, obtained by running time
backwards from $0$ to $-\infty$. By changing the support of the signal
spaces from the positive half-line to the negative half-line stability
defined through the forced response (BIBO stability) becomes
equivalent to asymptotic stability in the sense of Lyapunov for the
reversed time-direction as $t$ tends to {\mbox MINUS} infinity.
\begin{prop}\label{prop2} {\sf Let $\bf P$ be a linear system as in Proposition \ref{prop1} with $x(0)=0$. Then $y\in \cL_2(-\infty,0]$ for all $u\in \cL_2(-\infty,0]$ if and only if the matrix $-A$ is Hurwitz.}
\end{prop}

We now turn to the situation where inputs and outputs may have support
on the doubly-infinite time-axis. In this case the following holds, e.g.\ see \cite[p. 101]{zdg}.
\begin{prop}\label{prop3} {\sf Let $\bf P$ be a linear system as in Proposition \ref{prop1}. Then there exists $y\in \cL_2(-\infty,\infty)$ for all $u\in \cL_2(-\infty,\infty)$ if and only if $A$ has no imaginary-axis eigenvalues. Moreover, if this condition holds, $y$ is determined uniquely by $\hat{y}(s)=(C(sI-A)^{-1}B+D)\hat{u}(s)$.}
\end{prop}

We remark that Proposition \ref{prop3} is the natural generalisation
of Proposition \ref{prop1} when systems are viewed as operators.  A
linear system in Proposition \ref{prop1} becomes a multiplication
operator on the Fourier transformed spaces.  The operator is bounded
if and only if the ``symbol'' (the transfer-function) belongs to
$H_\infty$, which under the controllability and observability
assumption is equivalent to $A$ being Hurwitz.  On the double-axis a
multiplication operator on the Fourier transformed spaces is bounded
if and only if the symbol belongs to $L_\infty$---which for rational
symbols excludes only poles on the imaginary axis.

In Proposition \ref{prop3} there is no longer any relationship between
a notion of BIBO stability and Lyapunov stability (in either
time-direction). Clearly, both $A$ and $-A$ may fail to be Hurwitz.
Since only the existence of some $y\in\cL_2(-\infty,\infty)$ is
required for a given $u\in\cL_2(-\infty,\infty)$, and the free motion
solutions of (\ref{system}) are ignored, this is not surprising.
Propositions \ref{prop1} and \ref{prop2}, by contrast, establish a
connection between BIBO stability and Lyapunov stability as
$t\to+\infty$ (respectively, $t\to -\infty$) without putting in
explicit requirements on the free motion solutions.

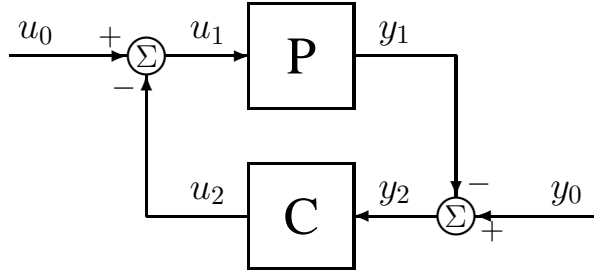
\begin{figure}
\begin{center}
\setlength{\unitlength}{.01in}%{0.012500in}
\parbox{304\unitlength}
{\begin{picture}(304,150)
\thicklines
\put(14,129){\makebox(0,0){{\large $u_0$}}}
\put(104,129){\makebox(0,0){{\large $u_1$}}}
\put(202,129){\makebox(0,0){{\large $y_1$}}}
\put(104,45){\makebox(0,0){{\large $u_2$}}}
\put(202,45){\makebox(0,0){{\large $y_2$}}}
\put(292,45){\makebox(0,0){{\large $y_0$}}}
\put(0,117){\vector(1,0){62}}
\put(53,127){\makebox(0,0){$+$}}
\put(72,117){\circle{20}}
\put(72,117){\makebox(0,0){$\Sigma$}}
\put(82,117){\vector(1,0){44}}
\put(180,117){\line(1,0){54}}
\put(126,90){\framebox(54,54){{\LARGE {\rm P}}}}
\put(234,117){\vector(0,-1){74}}
\put(246,50){\makebox(0,0){$-$}}
\put(72,33){\vector(0,1){74}}
\put(60,100){\makebox(0,0){$-$}}
\put(72,33){\line(1,0){54}}
\put(126,6){\framebox(54,54){{\LARGE {\rm C}}}}
\put(224,33){\vector(-1,0){44}}
\put(234,33){\circle{20}}
\put(234,33){\makebox(0,0){$\Sigma$}}
\put(253,25){\makebox(0,0){$+$}}
\put(306,33){\vector(-1,0){62}}
\end{picture}
}
\end{center}
\caption{Standard feedback configuration.}
\label{Figure1}
\end{figure}

We now consider the feedback interconnection in the form of Fig.\ 
\ref{Figure1} where $\bP$ and $\bC$ are linear systems.  The existence
of signals $u_i,y_j$ ($i,j\in\{1,2\}$) in $\cL_2[0,\infty)$ which
satisfy the feedback equations for a given pair of external inputs
$u_0,y_0$ in $\cL_2[0,\infty)$, for a given set of initial conditions, is
a well-known and natural definition of stability in terms of the
forced response.  From Proposition \ref{prop1} stability in this sense
is equivalent to asymptotic stability in the sense of Lyapunov of the
combined state-space (assuming minimal realizations for $\bP$ and
$\bC$ and well-posedness).  Again, BIBO stability inherits the
required time-asymmetry from the asymmetry of the support interval
$[0,\infty)$.

It is apparent that the corresponding definition of BIBO stability for
this feedback interconnection with $\cL_2(-\infty,\infty)$ signals,
generalising Proposition \ref{prop3}, will not correspond to a
sensible notion of closed-loop stability.  Indeed, we can easily check
that a system $\bP$ with transfer funtion $P(s)=1/(s-1)$ is
``stabilised'' by any of the controllers with $C(s)=2$, $C(s)=0$, or
$C(s)=-0.5/(s+1)$.  (In conventional terms the controllers give
closed-loop poles which are in the open left-half plane (LHP), the
open right-half plane (RHP), and in both half planes, respectively.)

We can summarize the points so far as follows.  Stability is a
time-asymmetric concept---the requirement of an asymptotic property as
$t$ tends to \mbox{PLUS} infinity defines a time arrow.  If stability
is defined by requiring bounded outputs in response to bounded inputs
then a time arrow is not obviously implied.  However, for signal
spaces with support on a positive (resp.\ negative) half-line, the
definition turns out to imply a positive (resp.\ negative) time arrow.
On the other hand, a bounded-input bounded-output definition of
stability for signals with support on the doubly-infinite time-axis
does not define a preferred time arrow. Stable systems defined by bounded ``multiplication operators''
may be stable in the sense of Lyapunov in the positive time-direction,
in the negative time-direction or in neither direction.

\subsection{The two-sided time axis and causality}\label{sectionIIB}
The fact that the doubly-infinite time axis causes problems for the
analysis of stability and of stabilisation was pointed out in
\cite{GS95}. The explanation given there is consistent with that of
Section \ref{section:timeasymmetry}, but the overall argument was
somewhat different. We now summarize the reasoning of \cite{GS95}.

Two systems ${\bf P}_i$ ($i=1,2$) defined by convolution operators
were considered:
\[
y(t)=\int_{-\infty}^\infty h_i(t-\tau)u(\tau)d\tau = h_i*u
\]
where $h_1(t)=e^t$ for $t\geq 0$ and zero otherwise, and $h_2(t)=-e^t$
for $t\leq 0$ and zero otherwise, respectively. Each system has
(double-sided Laplace) transfer function equal to $1/(s-1)$, but with differing
regions of convergence. The first system is unstable and causal and
the second is stable and non-causal (in fact anticausal) according to
the usual definitions.

When viewed on $\cL_2(-\infty,\infty)$, ${\bf P}_2$ is a bounded
operator and hence is a stable system in an input-output sense. On the
other hand, it was shown in \cite{GS95} that ${\bf P}_1$ fails to be
stabilisable on $\cL_2(-\infty,\infty)$.  This is a counterintuitive
result since ${\bf P}_1$ is stabilisable in the ordinary way on any
positive half-line.  The proof that ${\bf P}_1$ fails to be
stabilisable on the doubly-infinite time-axis reduces to the
observation that the graph of ${\bf P}_1$ fails to be closed.

It was also pointed out in \cite{GS95} that the closure of the graph
of ${\bf P}_1$ coincides with the graph of ${\bf P}_2$.  Once the
graph is closed there appears to be no problem with stabilisation.  But in
closing the graph ``anti-causal'' trajectories are brought in which
are inconsistent with the convolution representation of the system, so
this was considered inadmissible.

Another possible remedy discussed in \cite{GS95} was to consider the
underlying differential equation representations rather than the
convolution representations.  In fact both systems are defined by the
same differential equation
\begin{equation}\label{differentialequation}
\dot{y} = y + u.
\end{equation}
More precisely, the trajectories of both ${\bf P}_1$ and ${\bf P}_2$
satisfy this equation.  In terms of ``flow of time'' thinking,
$\bP_1$ appears to arise by solving this
equation forwards in time while $\bP_2$ is obtained by solving
it backwards. This suggestion seems to make stronger the argument to
consider ${\bf P}_1$ and ${\bf P}_2$ to be the same system.  But this
was considered unnatural in \cite{GS95} on the grounds
that it appears to abandon any notion of causality, or that it leaves
the direction of time undefined.

The discussion of Section \ref{section:timeasymmetry} allows the
difficulties pointed out in \cite{GS95} to be explained in a new way.
Let us suppose we are willing to accept the closure of the graph of
${\bf P}_1$ which makes it ``stabilizable'' on the double-axis in a bounded-input/bounded-output sense. As
explained, ${\bf P}_1$ and ${\bf P}_2$ can now be thought of as one
and the same system defined by (\ref{differentialequation})---a
state-space description as in (\ref{system}-\ref{system1}) solved forwards or
backwards as desired. Does the closure of the graph resolve the
difficulty pointed out in \cite{GS95}? The answer is no, since the notion of
stability does not correspond to the usual notions. As is made clear
by Proposition \ref{prop3}, the feedback system may turn out to be
stable in a conventional sense in the forward, backward or neither
time-directions.

\subsection{The work of M\"{a}kil\"{a}, Partington and Jacob}

A number of interesting observations and contributions have followed
from \cite{GS95} which we would like to comment on here.

The fact that a causal system on the double-axis can have a non-causal
closure has led to a study of ``closability'' and ``causal
closability'' as questions in their own right. M\"{a}kil\"{a}
\cite{makila9} has shown that the lack of causal closability for the
example of \cite{GS95} extends to general $\cL_p$ spaces on the
double-axis. Jacob and Partington \cite{jacob3} give general
characterisations of the graphs of time-invariant systems and derive
necessary and sufficient conditions for the closure of a closable
system to be causal. In \cite{makila1} M\"{a}kil\"{a} and Partington
consider weighted $\cL_2$-spaces on the double-axis and show that,
when signals have very rapid decrease to zero towards $-\infty$,
causal convolution operators may be closed operators. (So there is no
issue of causality being lost due to the operation of closure.)

On the question of stabilization on the double time-axis, Jacob
\cite{jacob} has made an interesting suggestion. We have seen already
that closing the graph and applying the BIBO stability definition
fails to recover the usual concept of stability. Jacob proposed that
causality of the closed-loop operators of the feedback system be added
as an extra requirement. Jacob showed that the resulting
characterisation of stability agrees with the usual definitions for
linear time-invariant systems. In the context of the present paper we
can re-interpret this result by saying that the causality condition
forces the positive time-arrow into feedback system stability. We can
understand this as follows. In \cite{jacob3} it is shown that a closed
linear time-invariant system is causal on $\cL_2(-\infty,\infty)$ if
and only if the corresponding transfer function belongs to a certain
Smirnov class. For finite dimensional systems this is equivalent to
the transfer function having no right half-plane poles. Thus, in
Proposition \ref{prop3}, if $\bf P$ is required to be causal, BIBO
stability agrees with Lyapunov stability with the positive time-arrow.

M\"{a}kil\"{a} and Partington in \cite{makila1} make an interesting
observation on the possible extension of Jacob's idea to the
time-varying situation. They consider a causal, convolution operator
derived from the underlying differential equation
\begin{equation}
\dot{y}(t)+a(t)y(t)=u(t)\label{timevarying}
\end{equation}
where $a(t)=-1$ for $t\leq 0$ and $a(t)=+1$ for $t>0$ and point out
that the closure of the $\cL_2(-\infty,\infty)$ graph of the
convolution system is not the graph of an operator. Essentially this
boils down to the fact that there are free motion solutions
$y(t)=ce^{-|t|}$, $u(t)=0$, where $c$ is a constant, which can be
approximated arbitrarily closely by elements of the graph. This raises the
question of whether the approach of Jacob can recover a theory of
stabilization which is consistent with the single-axis case. At the
same time it is pointed out that the system is stabilizable in a Lyapunov sense by the feedback
\begin{equation}\label{feedbacknegativetwo}
u(t)=-2y(t).
\end{equation}
In the present context this example highlights the care that is needed
in defining stability for time-varying systems, even in the
conventional sense. The open-loop system (\ref{timevarying})
is Lyapunov stable in the forward time-direction for any initial
condition specified at any time (either positive or negative), but not
uniformly so. Incidentally, the same is true for stability in the
backwards time-direction. With the feedback
(\ref{feedbacknegativetwo}) in force the system becomes uniformly
stable in the sense of Lyapunov in the forward time-direction and
unstable in reverse. In the perspective of the present paper, any
method to force agreement between BIBO stability on the double-axis
and conventional notions (such as requiring causality of the closed
loop operators) might be seen as tantamount to directly imposing the
desired time arrow within the stability definition.

In several papers (e.g., \cite{makila7,makila1,makila5})
M\"{a}kil\"{a} and Partington have advocated the use of a two-operator
model for systems on the doubly infinite time-axis in the form
$Ay=Bu$, where $A,B$ are causal, bounded operators, in contrast to a
single-operator model $y=Pu$, where $P$ is causal and possibly
unbounded. Closed-loop stability is defined as the existence of a
causal, bounded inverse of the feedback system operator mapping system
inputs to exogenous disturbances. Since this definition incorporates a
causality requirement on the closed-loop system there is evidently a
close relationship between this idea and the approach of Jacob.

\section{Time-Asymmetry and Optimal Regulation}\label{optimal control}

This section focusses on the time-asymmetry of the definition of
stability and its implications in the context of optimal regulation.
Firstly, a time-conjugation operator will be defined as well as the concepts of f-stability and b-stability. Then the finite-horizon quadratic regulator problem will be considered for a system running forwards in time and backwards in time, and it will be shown that the optimal cost is generally different. The infinite-horizon (asymptotic) regulator will also be considered in the same way. It will be shown that the optimal
cost can be expressed in terms of the two extremal solutions of the appropriate
algebraic Riccati equation. The result shows that ease of optimal regulation depends on the time-direction.

\subsection{The time-conjugation operator, f-stability and b-stability}\label{definitionofstability}
Let $\bP$ denote a dynamical system described by the state-space equations in (\ref{system}-\ref{system1}), initialized at time zero
and running forwards in time.  \newcommand{\fJ}{{\mathfrak{J}}} Let
$\fJ$ denote the operation on $\bP$ which corresponds to solving
(\ref{system}-\ref{system1}) backwards from $t=0$ followed by a flip
of the time axis (so the new system runs forward again). More specifically
we set $t_1=-t$, so that
\[
\frac{d}{dt}=-\frac{d}{dt_1},
\]
and then replace $t_1$ by $t$ which results in
\begin{eqnarray*}
-\dot{x}&=&Ax+Bu, \mbox{ with }x(0)=x_0\\
y&=&Cx+Du
\end{eqnarray*}
for the system $\fJ(\bP)$. The effect on the transfer function is as
follows: if $\bP$ has transfer function $P(s)$, then $\fJ(\bP)$ has
transfer function $P(-s)$.

Define the system $\bP$ to be {\em f-stable} if $A$ is Hurwitz, and
define $\bP$ to be {\em b-stable} if $\fJ(\bP)$ is f-stable, or
equivalently, if $-A$ is Hurwitz.  It is immediately obvious that a
linear time-invariant system of the type (\ref{system}-\ref{system1})
can never be both f-stable and b-stable.  Similarly, a controller
which makes (\ref{system}-\ref{system1}) f-stable cannot make it
b-stable as well.

\subsection{The finite-horizon linear quadratic regulator}

Let $\bf P$ be a linear time-invariant system which is controllable and observable and described by  (\ref{system}-\ref{system1}), as before, with $D=0$ and $x(0)=x_0$. Consider the problem of the regulation of $\bf P$ with criterion
\[
J=\int_0^T(y(t)^\prime Qy(t) + u(t)^\prime Ru(t))dt +x(T)^\prime H x(T).
\]
This has solution
\begin{equation}
u(t)=-R^{-1}B^\prime S(t)x(t)  \label{star1new}
\end{equation}
where
\begin{equation}  \label{star2new}
-\dot S(t)=S(t)A+A^\prime S(t)-S(t)BR^{-1}B^\prime S(t)+C^\prime QC
\end{equation}
and $S(T)=H$, with optimal cost
$J_{\rm f,T}=x_0^\prime S(0)x_0$  \cite{andmoore}.
When $\bP$ runs backwards in time from $x(0)=x_0$ with cost
\[
J=\int_{-T}^0(y(t)^\prime Qy(t) + u(t)^\prime Ru(t))dt +x(-T)^\prime H x(-T)
\]
we can check that the optimal control is still given by (\ref{star1new}) where $S(t)$ satisfies (\ref{star2new}) with $S(-T)=-H$, and that the optimal cost is $J_{b,T}=-x_0^\prime S(0)x_0$. It can be readily verified that $S(0)$ (forward case) is in general different from $-S(0)$ (backward case), and so the optimal cost is different in the two cases, e.g., if $A=B=C=Q=R=T=1$ and $H=10$, then $S(0)=2.5415$ in the forward case and $-S(0)=0.5495$ in the backward case.

\subsection{The infinite-horizon linear-quadratic regulator}

Again let $\bf P$ be a linear time-invariant system which is
controllable and observable and described by
(\ref{system}-\ref{system1}) with $D=0$ and $x(0) = x_0$.  It is
well-known \cite{andmoore} that
\begin{equation}
  J = \int_0^\infty (y(t)'Qy(t) + u(t)'Ru(t)) dt
\label{J_lq}
\end{equation}
has a minimum given by $J_{\rm f,\infty}=x_0'S_+x_0$ where $S_+$ is the unique
positive-definite solution to the algebraic Riccati equation
\begin{equation}
A'S+SA-SBR^{-1}B'S +C'QC=0.
\label{ric}
\end{equation}
It is also well-known that $S_+$ is the unique solution of (\ref{ric})
for which $A-BR^{-1}B'S$ has all its eigenvalues in the open LHP.  In
the language of the present paper we can say that $S_+$ is the unique
solution of (\ref{ric}) which makes the system
(\ref{system}-\ref{system1}) f-stable with the controller $u =
-R^{-1}B'Sx$.

What happens if we require the minimisation of
\begin{equation}
  J = \int_{-\infty}^0 (y(t)'Qy(t) + u(t)'Ru(t)) dt
\label{J_lq_back}
\end{equation}
for (\ref{system}-\ref{system1}) running backwards in time?  This is
the same as the conventional problem for the system $\fJ(\bP)$.  It is
easy to see that the minimum is given by $J_{\rm b,\infty}=-x_0'S_-x_0$ where
$S_-$ is the unique negative-definite solution to (\ref{ric}).  It is
also well-known that $S_-$ is the unique solution of (\ref{ric}) for
which $A-BR^{-1}B'S$ has all its eigenvalues in the open RHP
\cite{zdg}.  In the language of the present paper we can say that
$S_-$ is the unique solution of (\ref{ric}) which makes the system
(\ref{system}-\ref{system1}) b-stable with the controller $u =
-R^{-1}B'Sx$.

In general $J_{\rm f,\infty}=x_0'S_+x_0$ and $J_{\rm b,\infty}=-x_0'S_-x_0$ are different.
This shows that ``difficulty of control'' is time-asymmetric for the
standard linear-quadratic regulator on the infinite horizon.  The
difference can be significant, e.g.\ if $A=1$, $B=\epsilon$, $C=1$,
$Q=1$ and $R=1$ then $S_+ = 2/\epsilon^2+1/2+O(\epsilon^2)$ and $S_- =
-1/2+O(\epsilon^2)$ for $\epsilon$ small.

\section{Time-Asymmetry and Modelling Uncertainty}\label{section:robustcontrol}
In this section we look at the topology for uncertainty in
feedback control and how this is affected by the time arrow. We will
see that dynamical systems which are close in the usual sense, that a
common controller can stabilise them and give a similar closed-loop
behaviour, may not be close if time is reversed.

\subsection{The gap metric and robustness of stability}\label{section:gapmetric}

Zames and El-Sakkary \cite{ZamesElSakkary} introduced a metric on
dynamical systems for the purpose of assessing robustness.  This was
based on the gap metric used in functional analysis to study
invertibility of operators \cite{KreinKrasnoselski,SzN47}. Specifically,
systems are considered to be operators on $\cL_2[0,\infty)$ with a
graph which is a closed subspace of $\cL_2[0,\infty)$. Consider two
linear systems $\bP_i$ ($i=1,2$) with transfer functions
\[
P_i(s)={n_i(s)}\left({m_i(s)}\right)^{-1}
\]
where $n_i(s)$ and $m_i(s)$ are coprime polynomials or, more
generally, right-coprime polynomial matrices. Let
\[
\left(n_i(-s)\right)^Tn_i(s)+\left(m_i(-s)\right)^Tm_i(s)=\left(d_i(-s)\right)^Td_i(s)
\]
with $\det(d_i(s))$ a Hurwitz polynomial and $(\;)^T$ representing
matrix transpose---the existence of such a polynomial (matrix)
$d_i(s)$ is a standard result in the theory of canonical factorization
\cite{youla61}.  Then,
\[
\cG_{\bP_i,\cH_2}:=\left(\begin{matrix}{m_i(s)}({d_i(s)})^{-1}\\ \\
{n_i(s)}({d_i(s)})^{-1}\end{matrix}\right)\cH_2 :=G_i(s)\cH_2
\]
is (the Fourier transform of) the graph of $\bP_i$, for $i=1,2$. Thus,
the {\em graph symbol} $G_i(s)$ generates the graph of $\bP_i$ as its
range. Then the gap between $\bP_1$ and $\bP_2$ is defined to be
$\delta_{\cH_2}(\bP_1,\bP_2):=\|\Pi_{\cG_{\bP_1,\cH_2}}-\Pi_{\cG_{\bP_2,\cH_2}}\|$
where $\Pi_{\mathcal K}$ denotes orthogonal projection onto a closed
subspace $\mathcal K$.

Let the feedback configuration of Fig.\ \ref{Figure1} be denoted by $[\bP,\bC]$, where $\bP$ and $\bC$ are linear systems defined as operators on $\cL_2[0,\infty)$ which may possibly be unbounded.
Define
\[
\bH_{\bP,\bC}:=\left(\begin{matrix}\bI\\\bP\end{matrix}\right)
(\bI-\bP\bC)^{-1}
\left(\begin{matrix}\bI & -\bC\end{matrix}\right)
\]
to be the operator mapping
$\left(\begin{array}{cc}u_0^T& y_0^T\end{array}\right)^T$ to $\left(\begin{array}{cc}u_1^T&y_1^T\end{array}\right)^T$. The following are basic robustness results for gap metric uncertainty.\\
%$\left(\begin{array}{c}u_0\\y_0\end{array}\right)$ to $\left(\begin{array}{c}u_1\\y_1\end{array}\right)$.

\begin{prop}\label{margin} \cite{GeS90} {\sf Assume that
the  closed-loop system $[\bP,\bC]$ is f-stable. Then, $[\bP_1,\bC]$ is f-stable
for all $\bP_1$ such that $\delta_{\cH_2}(\bP,\bP_1)\leq b$ if and only if $b<b_{\bP,\bC}$ where
\[
b_{\bP,\bC}:=\|\bH_{\bP,\bC}\|_\infty^{-1}.
\]
}\end{prop}
\vspace*{.1in}

\begin{prop}\label{continuity} \cite{ZamesElSakkary} {\sf  Assume that the closed-loop system $[\bP,\bC]$ is f-stable. Then, the following are equivalent:
\begin{itemize}
\item[(i)] $\delta_{\cH_2}(\bP_n,\bP)\to 0$ as $n\to \infty$.
\item[(ii)] $\bH_{\bP_n,\bC}$ is f-stable for sufficiently large $n$ and
$\|\bH_{\bP_n,\bC}-\bH_{\bP,\bC}\|_\infty\to 0$ as $n\to\infty$.
\end{itemize}}\end{prop}

Proposition \ref{continuity} was the primary justification for the
claim in \cite{ZamesElSakkary} that the gap metric defines the
``correct'' topology for robustness of feedback systems. In the
present context, it can be seen that the choice of a signal space with
support on the positive half-line is essential in achieving an
appropriate topology. To emphasize the point, if $\cL_2[0,\infty)$
were replaced by $\cL_2(-\infty,0]$ then the above proposition would
hold with f-stability replaced by b-stability.

Let us consider the case where systems are defined on $\cL_2(-\infty,\infty)$. Then we define
\[\delta_{\cL_2}(\bP_1,\bP_2):=\|\Pi_{\cG_{\bP_1,\cL_2}}-\Pi_{\cG_{\bP_2,\cL_2}}\|\] where
\[
\cG_{\bP_i,\cL_2}:=G_i(s)\cL_2
\]
and $\cL_2:=\cL_2(-j\infty,j\infty)$. With this definition,
$\cG_{\bP,\cL_2}$ is always closed, but may contain ``non-causal''
input-output pairs (as pointed out in \cite{GS95}---see also Section
\ref{sectionIIB}). It is easy to construct examples to demonstrate
that convergence of $\delta_{\cL_2}(\bP_n,\bP)$ to zero does not allow
any closed-loop stability prediction, e.g., $[\bP,\bC]$ f-stable does
not imply $[\bP_n,\bC]$ f-stable for sufficiently large $n$.

In \cite{glenn} Vinnicombe introduced a new metric
$\delta_v(\cdot,\cdot)$ on dynamical systems which defines the same
topology as $\delta_{\cH_2}(\cdot,\cdot)$, and which satisfies the
following inequality:
\[
\delta_{\cL_2}(\cdot,\cdot)\leq\delta_v(\cdot,\cdot)\leq \delta_{\cH_2}(\cdot,\cdot).
\]
The v-gap between $\bP_1$ and $\bP_2$ is defined as follows:
\begin{equation}\label{nugap}
\delta_v(\bP_1,\bP_2):=\left\{ \begin{array}{l}\delta_{\cL_2}(\bP_1,\bP_2)\mbox{ if }\\
\hspace*{.7cm} \wno(\det(G_2(-s)^TG_1(s)))=0,\\
1 \mbox{ otherwise,}\end{array}\right.
\end{equation}
where $\wno(g(s))$ denotes the winding number about the origin of
$g(s)$, as $s$ traces the standard Nyquist D-contour \cite{glenn,glenn_book}. 
A simple expression for $\delta_{\cL_2}(\cdot,\cdot)$ can be obtained using left fractional representations---let $P_i(s)=(\tilde{m}_i(s))^{-1}\tilde{n}_i(s)$ be a left-coprime polynomial fraction,  $\tilde{d}_i$ the Hurwitz polynomial matrices which satisfy
\[
\tilde n_i(s)\left(\tilde n_i(-s)\right)^T+\tilde m_i(s)\left(\tilde m_i(-s)\right)^T=\tilde d_i(s)(\tilde d_i(-s))^T,
\]
and define
\[
\tilde{G}_i(s):=\left(\begin{matrix} -({\tilde d_i(s)})^{-1}{\tilde n_i(s)},\; {({\tilde d_i(s)})^{-1}\tilde m_i(s)}\end{matrix}\right)
\]
for $i=1,2$. The graph of $\bP_i$ is the kernel of multiplication by $\tilde G_i(s)$ (in the respective space of signals $\cH_2$ or $\cL_2$).
The $\cL_2$-gap can now be expressed as
\[\delta_{\cL_2}(\bP_1,\bP_2) :=\|\tilde{G}_2(s)G_1(s)\|_\infty.
\]

It turns out that Propositions \ref{margin} and \ref{continuity} both hold with $\delta_{\cH_2}$ replaced by
$\delta_v$ (see \cite{glenn}). Since $\delta_{\cL_2}=\delta_v$ when
\begin{equation}\label{wno}
\wno(\det(G_2(-s)^TG_1(s)))=0
\end{equation}
holds, this condition effectively imposes a positive time-arrow on the
double-axis graph which forces f-stability to be retained under small
perturbations in $\delta_v(\cdot,\cdot)$. This is illustrated by the
following result (which can be readily derived from \cite[Theorem 4.2]{glenn};
see also \cite{shankwitz}).

\begin{prop}{\sf Let $[\bP,\bC]$ be f-stable and suppose $\delta_{\cL_2}(\bP_n,\bP)\to 0$ as $n\to\infty$. Then $[\bP_n,\bC]$ is f-stable for all sufficiently large $n$ if and only if $\wno(\det(G_n(-s)^TG(s)))=0$ for all sufficiently large $n$.
}\end{prop}

\subsection{The effect of the time-arrow on gap distances}\label{section:graphsandvgap}

We define a forward and a backward v-gap as follows,
\begin{eqnarray*}
\delta_{v,f}(\bP_1,\bP_2)&:=&\delta_v(\bP_1,\bP_2)\\
\delta_{v,b}(\bP_1,\bP_2)&:=&\delta_v(\fJ(\bP_1),\fJ(\bP_2)).
\end{eqnarray*}
It is straightforward to see that
\[\delta_{\cL_2}(\bP_1,\bP_2) =\delta_{\cL_2}(\fJ(\bP_1),\fJ(\bP_2)),
\]
so any difference between $\delta_{v,f}(\bP_1,\bP_2)$ and
$\delta_{v,b}(\bP_1,\bP_2)$ lies in the winding number condition in
(\ref{nugap}).  Let us examine this more closely. Note that
\[\det (G_2(-s)^TG_1(s))=\frac{h(s)}{\det(d_2(-s))\det(d_1(s))}
 \]
where
\begin{equation}\label{eq:h}
  h(s):=\det(m_2(-s)^Tm_1(s)+n_2(-s)^Tn_1(s)).
\end{equation}
If $\delta_{\cL_2}(\bP_1,\bP_2)<1$ then it can be shown that
$\wno(\det(G_2(-s)^TG_1(s)))$ is well-defined \cite{glenn}, in which
case $h(s)$ admits a canonical factorization
\begin{equation}\label{eq:hfactor} h(s)=h_+(s)h_-(s)\end{equation}
where $h_+(s)$ and $h_-(-s)$ are Hurwitz polynomials.
Thus, $\wno(\det(G_2(-s)^TG_1(s)))=0$ if and only if
\[
\deg(h_+(s))=\deg(\det(d_1(s))),
\]
or equivalently
\[
\deg(h_-(s))=\deg(\det(d_2(s))),
\]
It can be shown that the degree of $\det(\hat d_i(s))$ is equal to the
McMillan degree of $\bP_i$ (e.g.\ using the uniqueness of normalised
coprime factors over $H_\infty$ up to a constant unitary
transformation and the corresponding state-space realisations
\cite{mey87,vid88,zdg}).  Determining the graph symbol for
$\fJ(\bP_i)$ requires a canonical factorization
\[
\left(n_i(s)\right)^Tn_i(-s)+\left(m_i(s)\right)^Tm_i(-s)=(\hat
d_i(-s))^T\hat d_i(s)
\]
with $\det(\hat d_i(s))$ a Hurwitz polynomial. Again it can be shown
that the degree of $\det(\hat d_i(s))$ is equal to the McMillan degree
of $\bP_i$.  The corresponding winding number condition in
$\delta_{v,b}(\bP_1,\bP_2)$ can now be expressed as
\[\wno(\det((\hat d_2(-s)^{-1}){h(-s)}(\hat d_1(s)^{-1})))=0
\]
which is equivalent to $\deg(h_-(s))$ being equal to the McMillan
degree of $\bP_1$. We therefore obtain the following result.

\begin{prop}\label{prop7}{\sf Let $P_i(s)$ ($i=1,2$) be the rational transfer functions of linear time-invariant dynamical systems as above, with McMillan degrees $\mu_i$, and with $h,h_+,h_-$ as in (\ref{eq:h}-\ref{eq:hfactor}). Assume that $\delta_{\cL_2}(\bP_1,\bP_2)<1$.
\begin{itemize}
\item[1)] The following are equivalent:
\begin{itemize}
\item[a)] $\delta_{v,f}(\bP_1,\bP_2)<1$,
\item[b)] $\deg(h_+(s))=\mu_1$,
\item[c)] $\deg(h_-(s))=\mu_2$.
\end{itemize}
\item[2)] The following are equivalent:
\begin{itemize}
\item[a)] $\delta_{v,b}(\bP_1,\bP_2)<1$,
\item[b)] $\deg(h_-(s))=\mu_1$,
\item[c)] $\deg(h_+(s))=\mu_2$.
\end{itemize}
\item[3)] The following are equivalent:
\begin{itemize}
\item[a)] $\delta_{v,f}(\bP_1,\bP_2)=\delta_{v,b}(\bP_1,\bP_2)<1$,
\item[b)] $\mu_1=\mu_2=\deg(h_+(s))=\deg(h_-(s))$.
\end{itemize}
\end{itemize}
}\end{prop}

In the above proposition, {\em 1)} expresses the zero winding number
condition in (\ref{nugap}) in an equivalent form, while {\em 2)} does
the same for $\delta_{v,b}(\bP_1,\bP_2)$. It is interesting that when
the two conditions are combined as in {\em 3)} the result is a very
stringent requirement which includes the necessity that $P_1(s)$ and
$P_2(s)$ have the same McMillan degree. This serves to highlight the
fact that ``unmodelled dynamics'' which may account for a small error
in $\delta_{v,f}$ (and which may be neglected in the design of a
robust controller) will inevitably account for a substantial error in
$\delta_{v,b}$.

%The above result has the following rather entertaining corollary:
%if both $P_i(s)$ ($i=1,2$) are functions of $s^2$ with different McMillan degrees, e.g., $P_1(s)=1$ and $P_2(s)=\frac{1}{s^2}$, then $\delta_{v,f}(\bP_1,\bP_2)=\delta_{v,b}(\bP_1,\bP_2)=1$.
\begin{example} Consider two systems with different McMillan degrees, e.g.\  $P_1(s)=1$, $P_2(s)=1/s$. It can be computed that
$\delta_{v,f}(\bP_1,\bP_2)={1/\sqrt{2}}$. Proposition \ref{prop7} then tells us immediately that $\delta_{v,b}(\bP_1,\bP_3)=1$ where $P_3(s)=-1/s$. Similarly, if $P_4(s)={1/s^2}$, then Proposition \ref{prop7} tells us that both
$
\delta_{v,f}(\bP_1,\bP_4)=
\delta_{v,b}(\bP_1,\bP_4)=1$ since $P_4(s)=P_4(-s)$.
\end{example}

\section{Time-asymmetry and robust control}\label{section:robustcontrol2}
This section addresses the implications of the time-asymmetry in the theory of robust
control. In particular, we will also see that
a system which is ``easy'' to control in one direction of time may be
far from easy to control in the opposite direction.

\subsection{Optimal robustness and  difficulty of control}\label{section:difficulty}

In \cite{McGlover} it was shown that $b_{\bP,\bC}$ could be maximised
over all stabilising $\bC$ and that this amounts to solving a Nehari
problem \cite{zdg}. This optimum value, which we denote by
\[
b_{\rm opt, f}(\bP),
\]
can be interpreted as a measure of ease/difficulty of control, where a
value near to $1$ means the plant is ``easy to control'' and a value
near $0$ means the plant is ``hard to control''.

With the understanding that $b_{\rm opt, f}(\bP)$ has the meaning of
``ease of control'' with respect to the forward time-arrow for stabilty,
it is interesting to define
\[
b_{\rm opt, b}(\bP):=
b_{\rm opt, f}(\fJ(\bP)),
\]
which represents ``ease of control'' with respect to the backwards
time-arrow. Our main purpose in defining $b_{\rm opt, b}(\bP)$ is to
highlight the influence of the time-arrow in feedback
regulation.

Let $\bf P$ be a controllable and observable system which is described by the state-space equations in (\ref{system}-\ref{system1}), as before.
Then, following \cite{McGlover,GeS90},
\[
b_{\rm opt,f}({\bf P})=\sqrt{1-\lambda_{\rm max}(Y_+X_+)}
\]
where $Y_+$ is the positive definite solution of the Riccati equation
\begin{eqnarray}\nonumber
A_0Y+YA_0^*-YCR^{-1}C^*Y\phantom{xxxxxxxx}\\
+B(I-D^*R^{-1}D^*)B^*=0\label{riccatieq}
\end{eqnarray}
where $A_0=A-BD^*R^{-1}C$ and $R=I+DD^*$,
and $X_+$ is the corresponding solution to the ($Y$-dependent) Lyapunov equation
\begin{eqnarray}\nonumber
(A_0-YC^*R^{-1}C)^*X+X(A_0-YC^*R^{-1}C)\\
+C^*R^{-1}C=0\label{lyapunoveq}
\end{eqnarray}
for $Y=Y_+$.
Similarly, it can be seen that
\begin{equation}\label{bopteq}
b_{\rm opt,b}({\bf P})=\sqrt{1-\lambda_{\rm max}(Y_-X_-)}
\end{equation}
where $Y_-$ is the negative definite solution of the Riccati equation
(\ref{riccatieq})
while $X_-$ is the corresponding solution to (\ref{lyapunoveq}) for $Y=Y_-$.

In the following two examples we will see situations where
$b_{\rm opt, f}(\cdot)$ and $b_{\rm opt, b}(\cdot)$ are very
different.

\begin{example}\label{nearcancellation} {\em Near pole-zero cancellations.}
  
  Consider $P(s)=1+\frac{\epsilon}{s+1}$. Letting $-A=C=D=1$ and $B=\epsilon$ in equations (\ref{riccatieq})-(\ref{bopteq}) gives
  \[
  b_{\rm opt,f}(\bP)=\sqrt{1-\frac{ \sqrt{1+\epsilon+\epsilon^2/2}-1-\epsilon/2}{2\sqrt{1+\epsilon+\epsilon^2/2}}}.
  \]
It follows that for small values of $\epsilon$,
\[
b_{\rm opt, f}(\bP)= 1-\frac{1}{32}\epsilon^2 + O(\epsilon^3)
\]
and hence,
$b_{\rm opt, f}(\bP)\to 1$ as $\epsilon\to 0$. On the other hand,
\[
  b_{\rm opt,b}(\bP)=\sqrt{1-\frac{ \sqrt{1+\epsilon+\epsilon^2/2}+1+\epsilon/2}{2\sqrt{1+\epsilon+\epsilon^2/2}}}.
  \]
which leads to
\[
b_{\rm opt, b}(\bP)= \frac{1}{4}|\epsilon| + O(\epsilon^2)
\]
for small values of $\epsilon$, and hence $b_{\rm opt,b}(\bP)\to 0$ as $\epsilon\to 0$.
This is accounted for by the
  fact that $P(s)$ has a near pole-zero cancellation in the LHP, which
  is innocuous for f-stabilisation, but highly challenging for
  b-stabilisation. The latter is equivalent to f-stabilisation of
  $P(-s)$, which has a troublesome near pole-zero cancellation in the
  RHP.
\end{example}

\begin{example}\label{bicycleexample}{\em Riding Bicycles.}

A feedback stability problem in everyday experience is bicycle
riding.  An elementary model to study rider-bicycle stability is given
in \cite{astrom} which gives the following transfer function from steering
angle input to tilt angle:
\begin{equation}
\alpha V \frac{s+\beta V}{s^2 -\gamma}
\label{bicycle_forward}
\end{equation}
where $\alpha,\beta,\gamma$ are positive constants and $V$ is the
forward speed. This model has one RHP pole, but the zero is in the
LHP.  As such, this plant is not too difficult to control.

Let us consider what happens if we try to ride the bicycle {\em
backwards in time}.  This corresponds to trying to stabilise the
plant $P(-s)$  {\em forwards in time}.  The model still has one RHP
pole, but the zero is also in the RHP, which makes stabilisation
much more difficult.  Indeed if $V\beta=\sqrt{\gamma}$ the plant is technically not
stabilisable.
It is interesting to note that an experimental bicycle with the
steered wheel at the rear instead of the front has a transfer function
from steering angle input to tilt angle given by \cite{astrom} (see also \cite{limebeer})
\begin{equation}
\alpha V \frac{-s+\beta V}{s^2 -\gamma}.
\label{bicycle_backward}
\end{equation}
This is exactly the transfer function for the conventional
bicycle ridden {\em backwards in time}.

Figure \ref{figure2} shows the value of $b_{\rm opt,f}$ and $b_{\rm opt,b}$ versus $V$ with parameter values $\alpha=1/3$, $\beta=2$ and $\gamma=9$ (which are deemed reasonably realistic).  Recall that $b_{\rm opt,b}$ is the same as $b_{\rm opt,f}$ for the rear-wheel steered bicycle model
(\ref{bicycle_backward}) at the same $V$.  It can be observed that
$b_{\rm opt,b}$ is less than $b_{\rm opt,f}$ for any $V$.  Also,
$b_{\rm opt,b}$ is very small for low $V$, indicating difficulty of
control, and zero at $V=1.5$ m/s.  For larger $V$, $b_{\rm opt,b}$
increases, indicating that control becomes easier.  These results are
equivalent to the rear-wheeled steered bicycle being more difficult to
ride than the front-wheel steered one, but still being reasonably
controllable at higher speeds \cite{astrom}.

\begin{figure}[htb]
\begin{center}\psfrag{V}[c]{{$V$ in m/s}}
\psfrag{B}[b]{$b_{\rm opt}$}
\includegraphics[width=7.4cm]{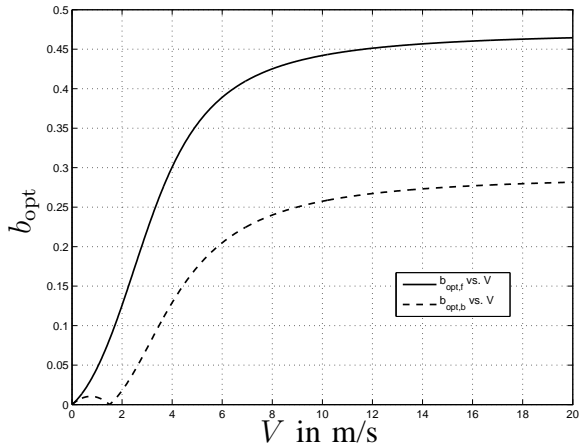}
\caption{$b_{\rm opt,f}$ and $b_{\rm opt,b}$ versus $V$ for the bicycle model of (\ref{bicycle_forward}) with $\alpha=1/3$, $\beta=2$, and $\gamma=9$.}
\label{figure2}
\end{center}
\end{figure}

\end{example}

\subsection{Time irreversible feedback phenomena}\label{geddanken}
The concept of ease or difficulty of control gives a thought-provoking
perspective on reversibility. Systems which in a limiting situation
are very difficult to control (in the sense that $b_{\rm opt,f}({\bf
  P})$ tends to zero) are unlikely to be observed in nature or
technology. Nevertheless, such a system may be easy to control in the
time-reversed direction (see Examples \ref{nearcancellation} and
\ref{bicycleexample}). This is independent of the fact that the
underlying differential equation can be integrated equally well in
either time-direction. This is reminiscent of phenomena (such as a
bottle falling from the table and shattering into many pieces) that
appear to be associated with an intrinsic direction of time even
though classical physics would also allow the reversed motion as a
solution (see Section \ref{physics} for a further discussion).

We expand this point in the context of Example \ref{bicycleexample}.
The loss of stabilizability of the rear-wheeled steered bicycle at
$V=\sqrt{\gamma/\beta}$ has the following interesting consequence.
Imagine a video of a rear-wheeled steered bicycle being ridden stably
at this critical speed.  Let us assume that it is possible to verify
from the video the actual speed (e.g., by knowing the frame-rate and
observing markings on the ground). An observer with a good grounding
in control theory would be led to the inescapable conclusion that the
video had been made when the said bicycle was actually being ridden
backwards in space (i.e., with a negative $V$) and then played
backwards in time as well, giving the impression of a forward motion.

\section{The arrow of time in physics}\label{physics}

The subject of the ``arrow of time'' is a well-known conundrum %??? and deep subject
in physics. The second law of thermodynamics states that the entropy
of a system increases with time. It is the time-asymmetry in this law
which gives rise to the notion of the ``thermodynamic arrow of time''.
The classical derivation of the second law in statistical mechanics
due to Boltzmann is connected with a famous puzzle known as
Loschmidt's paradox \cite{wikipedia}. This essentially points out that the laws of
mechanics used in the derivation of the second law are {\em
  time-symmetric} whereas the conclusion is not. Evidently the
time-asymmetry creeps in through the statistical assumptions. An
illuminating discussion of this issue is given in \cite{price}. Other
arrows of time have also been defined, for example (i) the
``psychological arrow''---the direction in which time passes as
perceived by a sentient being \cite{hawking,schulman}, (ii) the ``cosmological arrow''---the
direction of time in which the universe is expanding. Hawking
\cite{hawking} argues that the thermodynamic and psychological arrows
are always aligned with each other but these need not always be
aligned with the cosmological arrow (though they are at present).
%\addtolength{\textheight}{-5.6cm}  

In this paper we have described the time-asymmetry in the definition
of control systems stability as a time-arrow. In the theory of
dynamical systems there is also the notion of passivity, which again
defines a time-arrow. For electrical circuits the time-arrow of
passivity can be seen in the behaviour of the resistor, in contrast to
the inductor and capacitor which are time-symmetric in their
operation. If the electrical resistor were to operate backwards in
time one would observe a resistor gathering low-grade heat from the
environment and charging up a battery. This behaviour would be
recognised as a violation of the second law of thermodynamics (see
\cite[pages 260, 390-2]{prigogine}). In a similar way, an ideal linear
damper operating backwards in time extracts low-grade heat from the
environment to create mechanical work, in violation of the second law.
It seems that the arrow of time in passive systems or circuits
coincides with, or is the same as, the thermodynamic arrow.

How does the arrow of time for control system stability relate to
other time arrows? It is highly unlikely that a control engineer who
is designing a control system for a plant will give even a moment's
thought to the preferred time arrow for control. Without expressing the thought, the designer will seek
decaying free motion solutions in the direction in which time is
perceived to be passing. In this way the arrow of time for
control could be said to coincide with the psychological arrow. On the
other hand, in biological systems, active control is ubiquitous. It is
less obvious that, for example, homeostasis in a cell is aligned with
the psychological arrow. Here we will be content to raise the question of whether the stability arrow for control
systems in general can be directly related to the thermodynamic arrow, e.g.\ by considering information flow or the effect of internal
energy sources.

Finally, from a purely mathematical point of view, we observe that the arrow of
time for control systems stability appears identical with the arrow of
time for passivity. This supports the conclusion that the arrow of time for control
systems stability always coincides with the thermodynamic
(and psychological) arrow.

\section{Feedback loops and time delays}\label{delays}

Let $\bD_\tau \;:\; x(t) \mapsto x(t-\tau)$ denote a time delay
operator.  It seems superfluous to say that $\bD_\tau$ is physically realisable for $\tau>0$. Indeed the delay is a common feature of communication and control systems. For $\tau<0$, $\bD_\tau$ is the ideal predictor which is not believed to be physically realisable as a ``real-time'' device. At first sight this ``fact'' appears to be self-evident, but its subtlety is revealed on closer examination---indeed, a rigorous justification appears not to be available at present. An insightful discussion of the issue of ``causation'' and its connection with the arrow of time is given in Price \cite[Chapter 6]{price}. Price's suggestion that the asymmetry of causation ``is a projection of our own temporal asymmetry as agents in the world'' \cite[page 264]{price} is similar to the view expressed by Bertrand Russell: ``The law of causality, I believe, like much that passes muster among philosophers, is a relic of a bygone age, surviving, like the monarchy, only because it is erroneously supposed to do no harm''. This prevalence in physics and philosophy of an anthropocentric explanation of causation sits in opposition to the belief of the unrealizability of a ``prediction machine'' out of physical components and processes, and suggests that a deeper analysis of the question is needed.

In this paper we will not attempt to further debate the origin and explanation of causation.
In the next section we will simply highlight the striking difference in behaviour of feedback loops with small delays versus predictors and confirm the difference using the forward-time gap metric.

\subsection{Feedback stability, delays and predictors}

Consider a feedback system which consists of an
integrator in series with a time delay and negative unity feedback.
The governing equation is
\begin{equation}\label{integratorwithdelay}
\dot{x}(t)+x(t-\tau)=d(t)
\end{equation}
where $d(t)$ denotes an external disturbance.  We set
$d(t)\equiv 0$ and consider the totality of all free motion solutions
of the system equations. If all solutions decay as $t\to +\infty$ we say the system is f-stable. This definition agrees with the one given in Section
\ref{definitionofstability} for finite-dimensional systems.

For $\tau\geq 0$ we can verify that (\ref{integratorwithdelay}) is f-stable.
Taking Laplace transforms in (\ref{integratorwithdelay}) gives
\[\hat{x}(s)=\frac{1}{s+e^{-s\tau}}\hat{d}(s).
\]
We can verify that all zeros of $s+e^{-s\tau}=0$ are in the LHP so the
system is f-stable.

Now consider the case where $\tau<0$. Note that this corresponds to an
integrator with a predictor in negative feedback, which we would not
expect to be realizable in the forward time direction. In fact, $s+e^{-s\tau}$ has
infinitely many zeros in the RHP for any $\tau<0$ and hence the system
fails to be f-stable. It is evident that the system displays a discontinuity in the
asymptotic (as $t\to\infty$) behaviour of the free motion at the point
$\tau=0$.

Let us now consider the closeness of the systems involved using the
v-gap metric. Let $\bP$ denote the integrator and $\bP_\tau$ denote
the integrator in series with $\bD_\tau$. Regarding these as operators
on $\cL_2[0,\infty)$ we have the graph:
\[
{\cal G}_{\bP_\tau,\cH_2}=\left(\begin{matrix}\frac{s}{s+1}\\\frac{e^{-s\tau}}{s+1}\end{matrix}\right)\cH_2
\]
for $\tau\geq 0$. Then
\[
\delta_{\cL_2}(\bP,\bP_\tau)=\|\frac{s}{(s+1)^2}(1-e^{-s\tau})\|_\infty
\]
which tends to zero as $\tau\to 0$. Also
\[
G_2(-s)^TG_1(s)=\frac{-s^2+e^{s\tau}}{-s^2+1}
\]
so providing $|\tau|<\pi$ there are no crossings of the negative real
axis of this function when $s=j\omega$. Hence,
$\wno(G_2(-s)^TG_1(s))=0$, for $\tau$ sufficiently small. This implies
\[
\delta_{v,f}(\bP,\bP_\tau)=\delta_{\cL_2}(\bP,\bP_\tau)
\]
for $\tau\geq 0$ and sufficiently small and
$\delta_{v,f}(\bP,\bP_\tau)\to 0$ as $\tau\to 0$.

Now consider the case of $\bP_\tau$ with $\tau<0$. Again regarding
$\bP_\tau$ as an operator on $\cL_2[0,\infty)$ we have:
\[
{\cal G}_{\bP_\tau}=\left(\begin{matrix}\frac{se^{s\tau}}{s+1}\\\frac{1}{s+1}\end{matrix}\right)\cH_2
\]
and
$
\delta_{\cL_2}(\bP,\bP_\tau)=\|\frac{s}{(s+1)^2}(e^{s\tau}-1)\|_\infty
$ which tends to zero as $\tau\to 0$. Also,
\[
G_2(-s)^TG_1(s)=\frac{-s^2e^{-s\tau}+1}{-s^2+1},
\]
which behaves like $e^{-s\tau}$ for large $s$, so the winding number
of this function is not zero and $\delta_{v,f}(\bP,\bP_\tau)=1$ for
$\tau<0$.

The above analysis with the gap agrees with the earlier conclusion on
f-stability. For $\tau\geq 0$, f-stability was retained for
sufficiently small $\tau$, but lost for any $\tau<0$. Now we have seen
that, as long as $\tau\geq 0$, there is a small error in
$\delta_{v,f}$, but for any $\tau<0$, $\delta_{v,f}(\bP,\bP_\tau)=1$.

Finally, it is interesting to mention that the tolerance of feedback loops to small time-delays is guaranteed by a well-known sufficient condition
that the high-frequency loop-gain of the feedback loop is smaller than one (\cite{barman,wstability,Willems,Zames})---a condition routinely met in practice. It is easy to check that robustness to an arbitrarily small ``parasitic predictor'' in the loop would be guaranteed theoretically by the loop-gain being greater than one at arbitrarily high frequencies---a condition that appears impossible to achieve in a real feedback system.

\section{Synopsis}

\begin{enumerate}
\item Stability is a time-asymmetric concept. The requirement of an
  asymptotic property as $t$ tends to \mbox{PLUS} infinity defines a
  time arrow.
  
  \item A stability definition which requires bounded outputs in
  response to bounded inputs does not obviously imply a time arrow.
  For signal spaces with support on a positive (resp.\ negative)
  half-line, the definition turns out to imply a positive (resp.\ 
  negative) time arrow.
  
\item A bounded-input bounded-output definition of stability for
  signals with support on the doubly-infinite time-axis does not
  define a preferred time arrow. Stable systems defined by bounded multiplication operators may be
  stable in the sense of Lyapunov in the positive time direction, in
  the negative time direction or in neither direction.
  
\item The fact that the closure of the graph of an unstable causal
  system may coincide with the graph of a stable anti-causal system on
  the doubly-infinite time-axis need not be a fundamental obstacle in
  developing a usable control theory on the doubly-infinite
  time-axis.
  
\item Any method which modifies the BIBO definition of stability on
  the doubly-infinite time-axis to agree with conventional stability
  notions could be interpreted as the imposition of a positive
  time-arrow.
  
\item A time-conjugation operator on systems was defined as well as
  the concepts of f-stability and b-stability.
  
\item Both the finite-horizon and infinite-horizon quadratic regulators give a different optimal cost for a system running forwards in time and backwards in time. In the infinite horizon case the optimal cost can be expressed in terms of the two extremal solutions of the appropriate algebraic Riccati equation.
  
\item The role of the positive time arrow in the gap metric measure of
  uncertainty for dynamical systems was highlighted. The usual
  $\cH_2$-gap metric inherits the positive time arrow by virtue of
  systems being defined as operators on the positive half-line. The
  $\cL_2$-gap metric, which is well known to define an inappropriate
  topology for robust control, does not have a preferred
  time-direction due to the underlying operators being defined on the
  double-axis. The v-gap metric may be interpreted as the $\cL_2$-gap
  with an imposed time-arrow.
  
\item A time-conjugated v-gap metric was defined to measure closeness
  for robust b-stabilisation. It was seen that closeness of systems in
  the forward and backwards directions is a strong condition which
  includes the requirement of equal McMillan degrees.
  
\item It was seen that ease or difficulty of control as measured by optimal robustness in the gap metric is a property
  that depends on the time-arrow.
  
\item The situation of a plant which is easy to control in
  one time-direction but impossible to control in the other shows that
  irreversibility can be intimately related to control.
  
\item An engineering perspective of control suggests a close link
  between the control system stability arrow and the psychological
  arrow. Unified mathematical frameworks for passive circuits and
  feedback control suggest a close link between the control system
  stability arrow and the thermodynamic arrow. The question was raised whether the stability arrow for control systems can be directly related to the thermodynamic arrow.
  
  \item The issue of the non-realizability of the pure predictor as a ``real-time'' device and the connection with the arrow of time was highlighted as well as the difficulty of establishing non-realizability rigorously. The strongly contrasting behaviour of feedback loops in the presence of arbitrarily small time-delays or predictors was pointed out.
\end{enumerate}

\section{Acknowledgement}
We are grateful to Jan Willems for helpful comments on an earlier draft.

%\spacingset{1.7}

\end{document}